\documentclass[14pt]{article}
\usepackage{mathrsfs}
\usepackage{amsthm}
\usepackage{amsmath}
\usepackage{graphicx}
\usepackage{color}
\usepackage{amsfonts}
\newtheorem{theorem}{Theorem}[section]

\newtheorem{lemma}[theorem]{Lemma}

\newtheorem{corollary}[theorem]{Corollary}
\numberwithin{equation}{section}
\normalsize

\begin{document}
\title{\textbf{Contact processes with random vertex weights on oriented lattices}}

\author{Xiaofeng Xue \thanks{\textbf{E-mail}: xuexiaofeng@ucas.ac.cn \textbf{Address}: School of Mathematical Sciences, University of Chinese Academy of Sciences, Beijing 100049, China.}\\ University of Chinese Academy of Sciences}

\date{}
\maketitle

\noindent {\bf Abstract:}

In this paper we are concerned with contact processes with random vertex weights on oriented lattices. In our model, we assume that each vertex $x$ of $Z^d$ takes i. i. d. positive random value $\rho(x)$. Vertex $y$ infects vertex $x$ at rate proportional to $\rho(x)\rho(y)$ when and only when there is an oriented edge from $y$ to $x$. We give the definition of the critical value $\lambda_c$ of infection rate under the annealed measure and show  that $\lambda_c=[1+o(1)]/(dE\rho^2)$ as $d$ grows to infinity. Classic contact processes on oriented lattices and contact processes on clusters of oriented site percolation are two special cases of our model.

\noindent {\bf Keywords:}
Contact process, random vertex weights, oriented lattice, critical value.

\section{Introduction}\label{section one}
In this paper we are concerned with contact processes with random vertex weights on oriented lattices. For $d$-dimensional oriented lattice $Z^d$, there is an oriented edge from $x$ to $x+e_i$ for each $x\in Z^d$ and $1\leq i\leq d$, where
\[
e_i=(0,~\ldots,~0,\mathop 1\limits_{i \text{th}},0,~\ldots,~0).
\]
For $x,y\in Z^d$, we denote by $x\rightarrow y$ when $y-x\in \{e_i\}_{1\leq i\leq d}$. We denote by $O$ the origin of $Z^d$.

Let $\rho$ be a positive random variable such that $P(\rho>0)>0$ and $P(\rho\leq M)=1$ for some $M>0$. Let $\{\rho(x)\}_{x\in Z^d}$ be i. i. d. random variables such that $\rho(O)$ and $\rho$ have the same distribution. When $\{\rho(x)\}_{x\in Z^d}$ is given, the contact process with random vertex weights on oriented lattice $Z^d$ is a spin system with state space $\{0,1\}^{Z^d}$ and flip rates function given by
\begin{equation}\label{equ 1.1}
c(x,\eta)=
\begin{cases}
1 &\text{~if~}\eta(x)=1,\\
\lambda\sum\limits_{y:y\rightarrow x}\rho(x)\rho(y)\eta(y)&\text{~if~} \eta(x)=0
\end{cases}
\end{equation}
for each $(x,\eta)\in Z^d\times \{0,1\}^{Z^d}$, where $\lambda>0$ is a positive parameter called the infection rate. More details of the definition of spin systems can be find in Chapter 3 of \cite{LIG1985}.

Intuitively, this contact process describes the spread of an infection disease. Vertices in state $0$ are healthy and vertices in state $1$ are infected. An infected vertex waits for an exponential time with rate one to become healthy. An healthy vertex $x$ may be infected by an infected vertex $y$ when and only when there is an oriented edge from $y$ to $x$. The infection between $y$ and $x$ occurs at rate in proportional to $\rho(x)\rho(y)$.

Please note that the assumption $P(\rho<M)=1$ ensures the existence of our process according to the basic theory constructed in \cite{LIG1972}.

The contact processes with random vertex weights is introduced by Peterson in \cite{Pet2011} on finite complete graphs. He proves that the infection rate $\lambda$ has a critical value $\lambda_c=\frac{1}{E\rho^2}$ such that the disease survives for a long time with high probability when $\lambda>\lambda_c$ or dies out quickly with high probability when $\lambda<\lambda_c$.

Recently, contact processes in random environments or random graphs is a popular topic. In \cite{Chat2009}, Chatterjee and Durrett show that contact processes on random graphs with power law degree distributions have critical value $0$. This result disproves the guess in \cite{Pas2001a, Pas2001b} that the critical value is strictly positive according to a non-rigorous mean-field analysis. In \cite{Pet2011}, Peterson shows that contact processes with random vertex weights on complete graphs have critical value $\frac{1}{E\rho^2}$, which is consistent with the estimation given by the mean-field calculation. In \cite{Chen2009, Yao2012}, Yao and Chen shows that complete convergence theorem holds for contact processes in a random environment on $Z^d\times Z^+$. The random environment they set includes the bond percolation model as a special case.

In our model, if $\rho$ satisfies that $P(\rho=1)=1-P(\rho=0)=p$, then our model can be regarded as contact processes on clusters of oriented site percolation on $Z^d$. In \cite{Kesten1990}, Kesten shows that site percolation on $Z^d$ has critical probability $[1+o(1)]/2d$. We are inspired a lot by this result.

\section{Main result}\label{section two}
Before giving our main results, we introduce some notations. We assume that the random variables $\{\rho(x)\}_{x\in Z^d}$ are defined on a probability space $\big(\Omega, \mathcal{F}, P\big)$. We denote by $E$ the expectation operator with respect to $P$.

For any $\omega\in \Omega$, we denote by $P_\lambda^\omega$ the probability measure of our contact process on oriented lattice $Z^d$ with infection rate $\lambda$ and vertex weights $\{\rho(x,\omega)\}_{x\in Z^d}$. The probability measure $P_\lambda^\omega$ is called the quenched measure. We denote by $E_\lambda^\omega$ the expectation operator with respect to $P_\lambda^\omega$. We define
\[
P_{\lambda,d}(\cdot)=E\big[P_\lambda^\omega(\cdot)\big],
\]
which is called the annealed measure. We denote by $E_{\lambda,d}$ the expectation operator with respect to $P_{\lambda,d}$.

For any $t\geq 0$, we denote by $\eta_t$ the configuration of our process at the moment $t$. In this paper, we mainly deal with the case that all the vertices are infected at $t=0$. In later sections, if we need deal with the case that
\[
A=\{x:\eta_0(x)=1\}\neq Z^d,
\]
then we will point out the initial infected set $A$ and write $\eta_t$ as $\eta_t^A$. When $\eta_t$ is with no upper script, we refer to the case that
\[
\{x:\eta_0(x)=1\}=Z^d.
\]
According to basic coupling of spin systems, it is easy to see that
\[
P_{\lambda,d}(\eta_t(O)=1)\leq P_{\lambda,d}(\eta_s(O)=1)
\]
for $t>s$ and
\[
P_{\lambda_1,d}(\eta_t(O)=1)\leq P_{\lambda_2,d}(\eta_t(O)=1)
\]
for $\lambda_1<\lambda_2$. As a result, it is reasonable to define the following critical value of the infection rate.
\begin{equation}\label{equ 2.1}
\lambda_c(d)=\sup\big\{\lambda:\lim\limits_{t\rightarrow+\infty}P_{\lambda,d}(\eta_t(O)=1)=0\big\}.
\end{equation}
Please note that our process is symmetric under the annealed measure $P_{\lambda,d}$. So, $P_{\lambda,d}(\eta_t(x)=1)$ does not rely on the choice of $x$. As a result, when $\lambda<\lambda_c(d)$,
\[
\lim_{t\rightarrow+\infty}P_{\lambda,d}(\eta_t(x)=1\text{~for some~}x\in A)=0
\]
for any finite $A\subseteq Z^d$ and hence $\eta_t$ converges weakly to the configuration that all the vertices are healthy as $t$ grows to infinity.

Our main result is the following limit theorem of $\lambda_c(d)$.
\begin{theorem}\label{theorem main 2.1}
Assume that $P(\rho>0)>0$ and $P(\rho<M)=1$ for some $M>0$, then
\begin{equation}\label{equ 2.2 main}
\lim_{d\rightarrow+\infty}d\lambda_c(d)=\frac{1}{E\rho^2}.
\end{equation}
\end{theorem}

Theorem \ref{theorem main 2.1} shows that the critical value $\lambda_c(d)$ is approximately inversely proportional to the dimension $d$, the ratio of which is the reciprocal of the second moment of $\rho$.

When $\rho\equiv1$, our process is the classic contact process on oriented lattice. In this case, we write $\lambda_c(d)$ as $\lambda_d$. When $\rho$ satisfies that
\[
P(\rho=1)=1-P(\rho=0)=p
 \]
for some $p\in (0,1)$, our process is the contact process on clusters of oriented site percolation on $Z^d$. In this case, we write $\lambda_c(d)$ as $\lambda_c(d,\text{site},p)$. There are two direct corollaries of Theorem \ref{theorem main 2.1}.

\begin{corollary}\label{corollary 2.2}
\begin{equation}\label{equ 2.3}
\lim_{d\rightarrow+\infty}d\lambda_d=1.
\end{equation}
\end{corollary}

We say $y$ is $x$'s neighbor when $y\rightarrow x$, then Corollary \ref{corollary 2.2} shows that $\lambda_d$ is approximately to the reciprocal of the number of neighbors. In \cite{Grif1983, Hol1981, Pem1992}, Holley, Liggett, Griffeath and Pemantle show that this conclusion holds for contact processes on non-oriented lattices and regular trees. In \cite{Xue2014}, Xue shows that the same conclusion holds for threshold one contact processes on lattices and regular trees.

\begin{corollary}\label{corollary 2.3} For $p\in (0,1)$,
\begin{equation}
\lim_{d\rightarrow+\infty}dp\lambda_c(d,\text{site},p)=1.
\end{equation}
\end{corollary}

Corollary \ref{corollary 2.3} shows that $\lambda_c=[1+o(1)]/(dp)$ as $d$ grows to infinity for contact processes on clusters of oriented site percolation. In \cite{Xue2014b}, Xue claims that the same conclusion holds for contact process on clusters of oriented bond percolation on $Z^d$.

\quad

Please note that the critical value $\lambda_c(d)$ we define is under the annealed measure $P_{\lambda,d}$. We can also define critical value $\lambda_c(\omega)$ under the quenched measure such that
\[
\lambda_c(\omega)=\sup\big\{\lambda:\forall~x\in Z^d, \lim_{t\rightarrow+\infty}P_\lambda^\omega(\eta_t(x)=1)=0\big\}
\]
for any $\omega\in \Omega$. $\lambda_c(\omega)$ is a random variable. However, according to the ergodic theorem for i. i. d. random variables, it is easy to see that
\[
P(\omega: \lambda_c(\omega)=\lambda_c(d))=1.
\]
So we only need to deal with the critical value under the annealed measure.

The proof of Theorem \ref{theorem main 2.1} is divide into two sections. In Section \ref{section three}, we will prove that
\[
\liminf_{d\rightarrow+\infty}d\lambda_c(d)\geq\frac{1}{E\rho^2}.
\]
In Section \ref{section four}, we will prove that
\[
\limsup_{d\rightarrow+\infty}d\lambda_c(d)\leq\frac{1}{E\rho^2}.
\]

\section{Lower bound}\label{section three}
In this section we give a lower bound of $\lambda_c(d)$.
\begin{lemma}\label{lemma 3.1}
For each $d\geq 1$,
\[
\lambda_c(d)\geq \frac{1}{dE\rho^2}
\]
and hence
\[
\liminf_{d\rightarrow+\infty}d\lambda_c(d)\geq \frac{1}{E\rho^2}.
\]
\end{lemma}

\proof

We use $f_t$ to denote
\[
E_{\lambda,d}\big[\rho(O)\eta_t(O)\big]=E\big[\rho(O,\omega)P_\lambda^\omega(\eta_t(O)=1)\big].
\]
According to the flip rates function of $\eta_t$ given by \eqref{equ 1.1}, $\rho(O,\omega)$ and $P_\lambda^\omega(\eta_t(O)=1)$ are positive correlated. Therefore,
\begin{equation*}
f_t\geq E\rho(O,\omega)E[P_\lambda^\omega(\eta_t(O)=1)]=E\rho P_{\lambda,d}(\eta_t(O)=1).
\end{equation*}
Hence,
\begin{equation}\label{equ 3.1}
P_{\lambda,d}(\eta_t(O)=1)\leq \frac{f_t}{E\rho}.
\end{equation}
Please note that the assumption $P(\rho>0)>0$ ensures that $E\rho>0$.

According to Hille-Yosida Theorem and \eqref{equ 1.1},
\begin{align*}
\frac{d}{dt}P_\lambda^\omega(\eta_t(O)=1)=&-P_\lambda^\omega(\eta_t(O)=1)\\
&+\lambda\sum_{y:y\rightarrow O}\rho(O)\rho(y)P_{\lambda}^\omega(\eta_t(O)=0,\eta_t(y)=1)\\
\leq& -P_\lambda^\omega(\eta_t(O)=1)+\lambda\sum_{y:y\rightarrow O}\rho(O)\rho(y)P_{\lambda}^\omega(\eta_t(y)=1).
\end{align*}
Therefore,
\begin{equation}\label{equ 3.2}
\frac{d}{dt}f_t\leq -f_t+\lambda\sum_{y:y\rightarrow O}E\big[\rho^2(O)\rho(y)P_\lambda^\omega(\eta_t(y)=1)\big].
\end{equation}
For each $y$ such that $y\rightarrow O$, $\eta_t(y)$ is only influenced by the vertices from which there are oriented pathes to $y$. Therefore, $\rho(O)$ is independent of $\rho(y)P_{\lambda}^\omega(\eta_t(y)=1)$ and hence
\begin{equation}\label{equ 3.3}
E\big[\rho^2(O)\rho(y)P_\lambda^\omega(\eta_t(y)=1)\big]=E\rho^2P_{\lambda,d}(\eta_t(y)=1)=E\rho^2f_t.
\end{equation}
By \eqref{equ 3.2} and \eqref{equ 3.3},
\begin{equation}\label{equ 3.4}
\frac{d}{dt}f_t\leq (d\lambda E\rho^2-1)f_t.
\end{equation}
According to Greenwood inequality and \eqref{equ 3.4},
\[
f_t\leq f_0\exp\{(d\lambda E\rho^2-1)t\}
\]
and hence
\begin{equation}\label{equ 3.5}
\lim_{t\rightarrow+\infty}f_t=0
\end{equation}
when $\lambda<\frac{1}{dE\rho^2}$.

Lemma \ref{lemma 3.1} follows \eqref{equ 3.1} and \eqref{equ 3.5}.

\qed

\section{Upper bound}\label{section four}
In this section we will prove that $\liminf\limits_{d\rightarrow+\infty}d\lambda_c(d)\geq \frac{1}{E\rho^2}$.

First we define the contact process $\widehat{\eta}_t$ where disease spreads through the opposite direction of the oriented edges. For any $\omega\in \Omega$, The flip rates of $\widehat{\eta}_t$ with random vertex weights $\{\rho(x,\omega)\}_{x\in Z^d}$ is given by
\[
\widehat{c}(x,\eta)=
\begin{cases}
1 &\text{~if~} \eta(x)=1,\\
\lambda\sum\limits_{y:x\rightarrow y}\rho(x)\rho(y)\eta(y) & \text{~if~}\eta(x)=0.
\end{cases}
\]
Hence, for $\widehat{\eta}_t$, $y$ may infect $x$ when and only when there is an edge from $x$ to $y$.

According to the graphic representation of contact processes introduced by Harris in \cite{Har1978}. There is a dual of $\eta_t$ and $\widehat{\eta}_t$ such that
\begin{equation}\label{equ 4.1}
P_\lambda^\omega\big(\eta_t(O)=1\big)=P_\lambda^\omega\big(\widehat{\eta}_t^O\neq\emptyset\big),
\end{equation}
where $\widehat{\eta}_t^O$ is $\widehat{\eta}_t$ with that $\{x\in Z^d:\widehat\eta_0(x)=1\}=\{O\}$.

Please note that in \eqref{equ 4.1} we utilize the identification of $\widehat{\eta}_t^O$ with
\[
\big\{x\in Z^d:\widehat{\eta}_t^O(x)=1\big\}.
\]
Since $\{\rho(x)\}_{x\in Z^d}$ are i. i. d., the events $\{\eta_t^O\neq \emptyset\}$ and $\{\widehat{\eta_t}^O\neq \emptyset\}$ have the same distribution under the annealed measure $P_{\lambda,d}$. Therefore, according to \eqref{equ 4.1},
\begin{equation}\label{equ 4.2 dual}
P_{\lambda,d}\big(\eta_t(O)=1\big)=P_{\lambda,d}\big(\eta_t^O\neq \emptyset\big).
\end{equation}

To control the size of $\eta_t^O$ from below, we introduce a Markov process $\zeta_t$ with state space $\{-1,0,1\}^{Z^d}$. For given $\{\rho(x)\}_{x\in Z^d}$, $\zeta_t$ evolves as follows. For each $x\in Z^d$, if $\zeta(x)=-1$, then $x$ is frozen in the state $-1$ forever. If $\zeta(x)=1$, then the value of $x$ waits for an exponential time with rate one to become $-1$. If $\zeta(x)=0$, then the value of $x$ flips to $1$ at rate
\[
\lambda\sum_{y:y\rightarrow x}\rho(x)\rho(y)1_{\{\zeta(y)=1\}}.
\]

So for $\zeta_t$, when an infected vertex becomes healthy, then it is removed (in the state $-1$) and will never be infected.

We use $\zeta_t^O$ to denote $\zeta_t$ with $\{x\in Z^d:\zeta_0(x)=1\}=\{O\}$ and $\{x\in Z^d:\zeta_0(x)=-1\}=\emptyset$. According to the basic coupling of Markov processes, there is a coupling of $\eta_t$ and $\zeta_t$ under quenched measure $P_\lambda^\omega$ such that
\begin{equation}\label{equ 4.3}
\eta_t^O\supseteq \{x\in Z^d:\zeta_t^O(x)=1\}
\end{equation}
for any $t>0$.

We use $C_t$ to denote $\{x\in Z^d:\zeta_t^O(x)=1\}$. Then, by \eqref{equ 4.2 dual} and \eqref{equ 4.3},
\begin{equation}\label{equ 4.4}
\lim_{t\rightarrow+\infty}P_{\lambda,d}\big(\eta_t(O)=1\big)\geq P_{\lambda,d}\big(\forall~t, C_t\neq \emptyset\big).
\end{equation}

We give another description of $\{\forall~t,C_t\neq \emptyset\}$. Let $\{T_x\}_{x\in Z^d}$ be i. i. d. exponential times with rate $1$. For any $x\rightarrow y$, let $U_{xy}$ be exponential time with rate $\lambda \rho(x)\rho(y)$. We assume that all these exponential times are independent. For
\[
O=x_0\rightarrow x_1\rightarrow x_2\rightarrow\ldots\rightarrow x_n=x,
\]
if $U_{x_{j}x_{j+1}}\leq T_{x_j}$ for each $0\leq j\leq n-1$, then we say that there is an infected path with length $n$ from $O$ to $x$, which is denoted by $O\Rightarrow_n x$.

In the sense of coupling,
\[
\{O\Rightarrow_n x\}=\{\exists~t,x\in C_t\}.
\]
Let $I_n=\{x:O\Rightarrow_n x\}$ and $L_n$ be the set of infected pathes with length $n$ from $O$. $\{\forall~t, C_t\neq \emptyset\}$ is equivalent to that there are infinite many vertices which have ever been infected. Therefore,
\[
\{\forall~t,C_t\neq \emptyset\}=\{\forall~n, I_n\neq \emptyset\}.
\]
and
\begin{align}\label{equ 4.5}
P_{\lambda,d}\big(\forall~t,C_t\neq \emptyset\big)&=\lim_{n\rightarrow+\infty}P_{\lambda,d}\big(I_n\neq \emptyset\big)\notag\\
&=\lim_{n\rightarrow+\infty}P_{\lambda,d}\big(|L_n|>0\big)\notag\\
&\geq\limsup_{n\rightarrow+\infty}\frac{\big(E_{\lambda,d}|L_n|\big)^2}{E_{\lambda,d}|L_n|^2}
\end{align}
according to H\"{o}lder inequality.

To calculate $E_{\lambda,d}|L_n|$ and $E_{\lambda,d}|L_n|^2$, we utilize the simple random walk $S_n$ on oriented lattice $Z^d$ with $S_0=O$ and
\[
P(S_{n+1}-S_n=e_i)=\frac{1}{d}
\]
for $1\leq i\leq d$. Let $\{\widehat{S}_n\}_{n=0}^{+\infty}$ be an independent copy of $\{S_n\}_{n=0}^{+\infty}$. We assume that $\{S_n\}_{n=0}^{+\infty}$ and $\{\widehat{S}_n\}_{n=0}^{+\infty}$ are defined on probability space $(\widetilde{\Omega},\mathcal{G},\widetilde{P})$ and are independent of $\{\rho(x)\}_{x\in Z^d}$, $\{T_x\}_{x\in Z^d}$ and $\{U_{xy}\}_{x\rightarrow y}$. We denote by $\widetilde{E}$ the expectation operator with respect to $\widetilde{P}$.

For a given path $O\rightarrow x_1\rightarrow x_2\rightarrow\ldots\rightarrow x_n$,
\[
P_\lambda^\omega(U_{x_jx_{j+1}}<T_{x_j},\forall 0\leq j\leq n-1)=\prod_{j=0}^{n-1}\Big[\frac{\lambda\rho(x_j,\omega)\rho(x_{j+1},\omega)}{1+\lambda\rho(x_j,\omega)\rho(x_{j+1},\omega)}\Big].
\]
and hence
\begin{equation}\label{equ 4.7}
P_{\lambda,d}(U_{x_jx_{j+1}}<T_{x_j},\forall 0\leq j\leq n-1)=E\prod_{j=0}^{n-1}\Big[\frac{\lambda\rho(x_j)\rho(x_{j+1})}{1+\lambda\rho(x_j)\rho(x_{j+1})}\Big].
\end{equation}
It is obviously that the right hand side of \eqref{equ 4.7} does not rely on the oriented path $x_1,x_2,\ldots,x_n$ we choose.

As a result,
\begin{equation}\label{equ 4.6}
E_{\lambda,d}|L_n|=d^nE\prod_{j=0}^{n-1}\Big[\frac{\lambda\rho(S_j)\rho(S_{j+1})}{1+\lambda\rho(S_j)\rho(S_{j+1})}\Big]
\end{equation}
for any given first $n$ steps $(S_0,S_1,\ldots,S_n)$ of the simple random walk.

Please note that we write $E$ not $\widetilde{E}\times E$ in the right hand side of \eqref{equ 4.6}. We mean that the right hand side of \eqref{equ 4.6} is a random variable with respect to $\mathcal{G}$ and is a constant with probability one.

To calculate $E_{\lambda,d}|L_n|^2$, we introduce the following notations.
\begin{align*}
\tau_1&=\inf\{n\geq 0: S_n=\widehat{S}_n, S_{n+1}=\widehat{S}_{n+1}\},\\
\sigma_1&=\inf\{n>\tau_1:S_n=\widehat{S}_n,S_{n+1}\neq \widehat{S}_{n+1}\},\\
L_1&=\sigma_1-\tau_1+1,\\
\tau_2&=\inf\{n>\sigma_1: S_n=\widehat{S}_n, S_{n+1}=\widehat{S}_{n+1}\},\\
\sigma_2&=\inf\{n> \tau_2:S_n=\widehat{S}_n,S_{n+1}\neq \widehat{S}_{n+1}\},\\
L_2&=\sigma_2-\tau_2+1,\\
&\ldots\ldots\\
\tau_{k}&=\inf\{n>\sigma_{k-1}: S_n=\widehat{S}_n, S_{n+1}=\widehat{S}_{n+1}\},\\
\sigma_{k}&=\inf\{n>\tau_k:S_n=\widehat{S}_n,S_{n+1}\neq \widehat{S}_{n+1}\},\\
L_k&=\sigma_k-\tau_k+1,\\
&\ldots\ldots\\
T&=\sup\{k:\tau_k<+\infty\}.
\end{align*}
Please note that $P(T<+\infty)=1$ for $d\geq 4$ according to the conclusion proven in \cite{Cox1983} that $P(\exists~n>0, S_n=\widehat{S}_n)<1$ for $d\geq 4$. Therefore, $\tau_k,\sigma_k,L_k$ are finite for $k\leq T$.

Furthermore, we define
\begin{align*}
A_0&=\{0\leq n<\tau_1:S_n=\widehat{S}_n,S_{n+1}\neq \widehat{S}_{n+1}\},\\
A_1&=\{\sigma_1<n<\tau_2:S_n=\widehat{S}_n,S_{n+1}\neq \widehat{S}_{n+1}\},\\
&\ldots\ldots\\
A_{T-1}&=\{\sigma_{T-1}<n<\tau_T:S_n=\widehat{S}_n,S_{n+1}\neq \widehat{S}_{n+1}\},\\
A_T&=\{n>\sigma_T:S_n=\widehat{S}_n,S_{n+1}\neq \widehat{S}_{n+1}\}.
\end{align*}
For $0\leq i\leq T$, we use $K_i$ to denote $|A_i|$.

After all this prepare work, we give a lemma which is crucial for us to give upper bound of $\lambda_c(d)$.
\begin{lemma}\label{Lemma 4.1}
Assume that $P(\rho>0)>0$ and $P(\rho<M)=1$. If $\lambda$ makes
\begin{equation}\label{equ 4.8}
\widetilde{E}\Bigg[\frac{2^{T+\sum\limits_{j=0}^TK_j}M^{6T+4\sum\limits_{j=0}^TK_j}\big(1+\lambda M^2\big)^{2\sum\limits_{j=1}^TL_j+2\sum\limits_{j=0}^TK_j}}{\lambda^{\sum\limits_{j=1}^TL_j-T}\big(E\rho^2\big)^{\sum\limits_{j=1}^TL_j+2T
+2\sum\limits_{j=0}^TK_j}}\Bigg]<+\infty,
\end{equation}
then
\[
\lambda_c(d)\leq \lambda.
\]
\end{lemma}
\proof For each $x\rightarrow z_1$ and $y\rightarrow z_2$, we define $F(x,y;z_1,z_2)$ as
\[
P_\lambda^\omega\big(U_{xz_1}\leq T_x,U_{yz_2}\leq T_y\big).
\]
By direct calculation,
\begin{equation}\label{equ 4.9}
F(x,y;z_1,z_2)
\begin{cases}
=\frac{\lambda^2\rho(x)\rho(y)\rho(z_1)\rho(z_2)}{[1+\lambda\rho(x)\rho(z_1)][1+\lambda\rho(y)\rho(z_2)]} &\text{~if~}x\neq y\text{~and~}z_1\neq z_2,\\
=\frac{\lambda^2\rho(x)\rho(y)\rho^2(z_1)}{[1+\lambda\rho(x)\rho(z_1)][1+\lambda\rho(y)\rho(z_1)]}&\text{~if~}x\neq y\text{~and~}z_1=z_2,\\
=\frac{\lambda\rho(x)\rho(z_1)}{1+\lambda\rho(x)\rho(z_1)}&\text{~if~}x=y\text{~and~}z_1=z_2,\\
\leq \frac{2\lambda^2\rho^2(x)\rho(z_1)\rho(z_2)}{[1+\lambda\rho(x)\rho(z_1)][1+\lambda\rho(x)\rho(z_2)]}&\text{~if~}x=y \text{~and~}z_1\neq z_2.
\end{cases}
\end{equation}
We denote by $P_n$ the set of all the oriented paths from $O$ with length $n$, then
\begin{align*}
E_{\lambda,d}|L_n|^2&=\sum_{\textbf{x}\in P_n}\sum_{\textbf{y}\in P_n}P_{\lambda,d}\big(\forall~0\leq i\leq n-1,U_{x_ix_{i+1}}\leq T_{x_i},U_{y_iy_{i+1}}\leq T_{y_i}\big)\\
&=\sum_{\textbf{x}\in P_n}\sum_{\textbf{y}\in P_n}EP_{\lambda}^{\omega}\big(\forall~0\leq i\leq n-1,U_{x_ix_{i+1}}\leq T_{x_i},U_{y_iy_{i+1}}\leq T_{y_i}\big)\\
&=\sum_{\textbf{x}\in P_n}\sum_{\textbf{y}\in P_n}E\big[\prod_{i=0}^{n-1}F(x_i,y_i;x_{i+1},y_{i+1})\big]\\
&=d^{2n}\sum_{\textbf{x}\in P_n}\sum_{\textbf{y}\in P_n}\frac{1}{d^{2n}}E\big[\prod_{i=0}^{n-1}F(x_i,y_i;x_{i+1},y_{i+1})\big]\\
&=d^{2n}(\widetilde{E}\times E)[\prod_{i=0}^{n-1}F(S_i,\widehat{S}_i;S_{i+1},\widehat{S}_{i+1})].
\end{align*}

Therefore, by \eqref{equ 4.6},
\begin{equation}\label{equ 4.10}
\frac{(E_{\lambda,d}|L_n|)^2}{E_{\lambda,d}|L_n|^2}=\Bigg\{\widetilde{E}\Bigg[\frac{E\prod\limits_{i=0}^{n-1}F(S_i,\widehat{S}_i;S_{i+1},\widehat{S}_{i+1})}{\Big(E\prod\limits_{i=0}^{n-1}
\frac{\lambda\rho(S_i)\rho(S_{i+1})}{1+\lambda\rho(S_i)\rho(S_{i+1})}\Big)^2}\Bigg]\Bigg\}^{-1}.
\end{equation}
Then by \eqref{equ 4.4} and \eqref{equ 4.5}, $\lambda\geq \lambda_c(d)$ when
\[
\limsup_{n\rightarrow+\infty} \widetilde{E}\Bigg[\frac{E\prod\limits_{i=0}^{n-1}F(S_i,\widehat{S}_i;S_{i+1},\widehat{S}_{i+1})}{\Big(E\prod\limits_{i=0}^{n-1}
\frac{\lambda\rho(S_i)\rho(S_{i+1})}{1+\lambda\rho(S_i)\rho(S_{i+1})}\Big)^2}\Bigg]< +\infty.
\]
Now we control
\[
\widetilde{E}\Bigg[\frac{E\prod\limits_{i=0}^{n-1}F(S_i,\widehat{S}_i;S_{i+1},\widehat{S}_{i+1})}{\Big(E\prod\limits_{i=0}^{n-1}
\frac{\lambda\rho(S_i)\rho(S_{i+1})}{1+\lambda\rho(S_i)\rho(S_{i+1})}\Big)^2}\Bigg]
\]
from above.

For the denominator $\Big(E\prod\limits_{i=0}^{n-1}
\frac{\lambda\rho(S_i)\rho(S_{i+1})}{1+\lambda\rho(S_i)\rho(S_{i+1})}\Big)^2$, if $S_i=\widehat{S}_i$ or $S_{i+1}=\widehat{S}_{i+1}$, then we narrow the factor $\frac{1}{1+\lambda\rho(S_i)\rho(S_{i+1})}$ to $\frac{1}{1+\lambda M^2}$, where $P(\rho<M)=1$ as we assumed.

For the numerator $E\prod\limits_{i=0}^{n-1}F(S_i,\widehat{S}_i;S_{i+1},\widehat{S}_{i+1})$, if $S_i=\widehat{S}_i$ or $S_{i+1}=\widehat{S}_{i+1}$, then we enlarge the factor $\frac{1}{1+\lambda\rho(S_i)\rho(S_{i+1})}$ to $1$. If $i\in A_k$ for some $k$, then by \eqref{equ 4.9}, we enlarge the factors
\[2\lambda^2\rho^2(S_i)\rho(S_{i+1})\rho(\widehat{S}_{i+1})\text{~and~} \lambda^2\rho(S_{i-1})\rho(\widehat{S}_{i-1})\rho^2(S_i)\]
to
\[2\lambda^2M^2\rho(S_{i+1})\rho(\widehat{S}_{i+1})\text{~and~} \lambda^2\rho(S_{i-1})\rho(\widehat{S}_{i-1})M^2.\]
If $i=\tau_k$ for some $k$, then we enlarge the factors
\[\lambda^2\rho(S_{i-1})\rho(\widehat{S}_{i-1})\rho^2(S_i)\text{~and~}\lambda\rho(S_i)\rho(S_{i+1})\]
to
\[\lambda^2\rho(S_{i-1})\rho(\widehat{S}_{i-1})M^2\text{~and~}\lambda M\rho(S_{i+1}).\]
If $i=\sigma_k$ for some $k$, then we enlarge the factors
\[
\lambda\rho(S_{i-1})\rho(S_i) \text{~and~} 2\lambda^2\rho^2(S_i)\rho(S_{i+1})\rho(\widehat{S}_{i+1})
\]
to
\[
\lambda\rho(S_{i-1})M \text{~and~} 2\lambda^2M^2\rho(S_{i+1})\rho(\widehat{S}_{i+1}).
\]
After all these operations, we can cancel many common factors in the numerator and denominator. For example, if $i,j\in A_k$ and $l\not\in A_k$ for each $i<l<j$, then we can abstract
\[
\Big[E\frac{\prod\limits_{l=1}^{j-i-1}\rho_l^2}{\prod\limits_{l=1}^{j-i-2}\big(1+\lambda\rho_l\rho_{l+1}\big)}\Big]^2
\]
from both numerator and denominator and cancel this common factor, where $\{\rho_l\}_{l=1}^{i-j-1}$ are i. i. d. and have the same distribution as that of $\rho$.

Therefore, after all the above operations, it is easy to see that
\begin{align*}
&\limsup_{n\rightarrow+\infty}\frac{E\prod\limits_{i=0}^{n-1}F(S_i,\widehat{S}_i;S_{i+1},\widehat{S}_{i+1})}{\Big(E\prod\limits_{i=0}^{n-1}
\frac{\lambda\rho(S_i)\rho(S_{i+1})}{1+\lambda\rho(S_i)\rho(S_{i+1})}\Big)^2}\\
&\leq \frac{2^{T+\sum\limits_{j=0}^TK_j}M^{6T+4\sum\limits_{j=0}^TK_j}\big(1+\lambda M^2\big)^{2\sum\limits_{j=1}^TL_j+2\sum\limits_{j=0}^TK_j}}{\lambda^{\sum\limits_{j=1}^TL_j-T}\big(E\rho^2\big)^{\sum\limits_{j=1}^TL_j+2T
+2\sum\limits_{j=0}^TK_j}}
\end{align*}
and
\begin{align}\label{equ 4.11}
&\limsup_{n\rightarrow+\infty} \widetilde{E}\Bigg[\frac{E\prod\limits_{i=0}^{n-1}F(S_i,\widehat{S}_i;S_{i+1},\widehat{S}_{i+1})}{\Big(E\prod\limits_{i=0}^{n-1}
\frac{\lambda\rho(S_i)\rho(S_{i+1})}{1+\lambda\rho(S_i)\rho(S_{i+1})}\Big)^2}\Bigg] \notag\\ &\leq\widetilde{E}\Bigg[\frac{2^{T+\sum\limits_{j=0}^TK_j}M^{6T+4\sum\limits_{j=0}^TK_j}\big(1+\lambda M^2\big)^{2\sum\limits_{j=1}^TL_j+2\sum\limits_{j=0}^TK_j}}{\lambda^{\sum\limits_{j=1}^TL_j-T}\big(E\rho^2\big)^{\sum\limits_{j=1}^TL_j+2T
+2\sum\limits_{j=0}^TK_j}}\Bigg].
\end{align}
Lemma \ref{Lemma 4.1} follows \eqref{equ 4.4}, \eqref{equ 4.5}, \eqref{equ 4.10} and \eqref{equ 4.11}.

\qed

\quad

Finally, we give the proof of $\limsup_{n\rightarrow+\infty}d\lambda_c(d)\leq \frac{1}{E\rho^2}$.

\proof[Proof of $~\limsup_{n\rightarrow+\infty}d\lambda_c(d)\leq \frac{1}{E\rho^2}$]

Let
\[
\tau=\inf\{n>0:S_n=\widehat{S}_n\}.
\]
Then according to (2.9) of \cite{Cox1983},
\[
P(2\leq \tau<+\infty)\leq \frac{C_1}{d^2},
\]
where $C_1$ does not depend on $d$. Therefore, according to strong Markov property,
\begin{align}\label{equ 4.12}
&P(T=m,K_i=k_i\text{~for~}0\leq i\leq m,L_i=l_i\text{~for~}1\leq i\leq m)\notag\\
&\leq\big(\frac{C_1}{d^2}\big)^{\sum\limits_{i=0}^mk_i+m-1}\big(\frac{1}{d}\big)^{\sum\limits_{i=1}^ml_i-m}
\end{align}
for all possible $m,k_i,l_i$. Please note that $k_0$ may take $0$ but $l_i\geq 1$ and $k_i\geq 1$ for $1\leq i\leq m$.

Let
\[
\lambda=\frac{\gamma}{dE\rho^2}
\]
for fixed $\gamma>1$, then by \eqref{equ 4.12},
\begin{align}\label{equ 4.13}
&\widetilde{E}\Bigg[\frac{2^{T+\sum\limits_{j=0}^TK_j}M^{6T+4\sum\limits_{j=0}^TK_j}\big(1+\lambda M^2\big)^{2\sum\limits_{j=1}^TL_j+2\sum\limits_{j=0}^TK_j}}{\lambda^{\sum\limits_{j=1}^TL_j-T}\big(E\rho^2\big)^{\sum\limits_{j=1}^TL_j+2T
+2\sum\limits_{j=0}^TK_j}}\Bigg]\notag\\
&\leq\sum_{m=0}^{+\infty}\sum_{k_0=0}^{+\infty}\sum_{k_1=1}^{+\infty}\cdots\sum_{k_m=1}^{+\infty}\sum_{l_1=1}^{+\infty}\cdots\sum_{l_m=1}^{+\infty}
\big(\frac{C_1}{d^2}\big)^{\sum\limits_{i=0}^mk_i+m-1}
\big(\frac{1}{d}\big)^{\sum\limits_{i=1}^ml_i-m}\notag\notag\\
&\times\frac{2^{m+\sum\limits_{j=0}^mk_j}M^{6m+4\sum\limits_{j=0}^mk_j}\big(1+\lambda M^2\big)^{2\sum\limits_{j=1}^ml_j+2\sum\limits_{j=0}^mk_j}}{\lambda^{\sum\limits_{j=1}^ml_j-m}\big(E\rho^2\big)^{\sum\limits_{j=1}^ml_j+2m
+2\sum\limits_{j=0}^mk_j}}.\notag\\
&=\sum_{m=0}^{+\infty}\sum_{k_0=0}^{+\infty}\Big(\frac{2C_1M^6\lambda}{d(E\rho^2)^2}\Big)^m\Big[\frac{2C_1M^4(1+\lambda M^2)^2}{d^2(E\rho^2)^2}\Big]^{k_0}\\
&\times\Big[\sum_{l=1}^{+\infty}\big(\frac{2C_1M^4(1+\lambda M^2)^2}{d^2(E\rho^2)^2}\big)^l\Big]^m\Big[\sum_{l=1}^{+\infty}\big(\frac{(1+\lambda M^2)^2}{d\lambda E\rho^2}\big)^l\Big]^m\frac{d^2}{C_1}.\notag
\end{align}

Since $\lambda=\frac{\gamma}{dE\rho^2}$ for $\gamma>1$,
\[
\frac{2C_1M^6\lambda}{d(E\rho^2)^2}\leq \frac{C_2}{d^2}
\]
and
\[
\frac{2C_1M^4(1+\lambda M^2)^2}{d^2(E\rho^2)^2}\leq \frac{C_3}{d^2}
\]
for sufficiently large $d$, where $C_2$ and $C_3$ does not depend on $d$ (but may depend on $\gamma$ and $\rho$).

We choose $\widehat{\gamma}$ such that $1<\widehat{\gamma}<\gamma$, then for sufficiently large $d$,
\[
\frac{(1+\lambda M^2)^2}{d\lambda E\rho^2}=\frac{\big(1+\frac{\gamma M^2}{dE\rho^2}\big)^2}{\gamma}<\frac{1}{\widehat{\gamma}}.
\]

Then, by \eqref{equ 4.13},
\begin{align}\label{equ 4.14}
&\widetilde{E}\Bigg[\frac{2^{T+\sum\limits_{j=0}^TK_j}M^{6T+4\sum\limits_{j=0}^TK_j}\big(1+\lambda M^2\big)^{2\sum\limits_{j=1}^TL_j+2\sum\limits_{j=0}^TK_j}}{\lambda^{\sum\limits_{j=1}^TL_j-T}\big(E\rho^2\big)^{\sum\limits_{j=1}^TL_j+2T
+2\sum\limits_{j=0}^TK_j}}\Bigg] \notag\\
&\leq \frac{d^2}{C_1}\sum_{m=0}^{+\infty}\sum_{k_0=0}^{+\infty}\big(\frac{C_2}{d^2}\big)^m\big(\frac{C_3}{d^2}\big)^{k_0}
\big[\sum_{l=1}^{+\infty}(\frac{C_3}{d^2})^l\big]^m\big[\sum_{l=1}^{+\infty}(\frac{1}{\widehat{\gamma}})^l\big]^m\notag.
\end{align}
For sufficiently large $d$,
\[
\sum_{l=1}^{+\infty}\big(\frac{C_3}{d^2}\big)^l=\frac{C_3}{d^2-C_3}\leq \frac{C_4}{d^2}
\]
and
\[
\sum_{k_0=0}^{+\infty}\big(\frac{C_3}{d^2}\big)^{k_0}=\frac{d^2}{d^2-C_3}\leq 2,
\]
where $C_4$ does not depend on $d$.

Therefore,
\begin{align}\label{equ 4.15}
&\widetilde{E}\Bigg[\frac{2^{T+\sum\limits_{j=0}^TK_j}M^{6T+4\sum\limits_{j=0}^TK_j}\big(1+\lambda M^2\big)^{2\sum\limits_{j=1}^TL_j+2\sum\limits_{j=0}^TK_j}}{\lambda^{\sum\limits_{j=1}^TL_j-T}\big(E\rho^2\big)^{\sum\limits_{j=1}^TL_j+2T
+2\sum\limits_{j=0}^TK_j}}\Bigg] \notag\\
&\leq \frac{2d^2}{C_1}\sum_{m=0}^{+\infty}\big(\frac{C_2}{d^2}\big)^m\big(\frac{C_4}{d^2}\big)^m\big[\frac{1}{\widehat{\gamma}-1}\big]^m\notag\\
&=\frac{2d^2}{C_1}\sum_{m=0}^{+\infty}\big[\frac{C_2C_4}{d^4(\widehat{\gamma}-1)}\big]^m.
\end{align}

For sufficiently large $d$, $\frac{C_2C_4}{d^4(\widehat{\gamma}-1)}<1$ and therefore
\begin{equation}\label{equ 4.16}
\widetilde{E}\Bigg[\frac{2^{T+\sum\limits_{j=0}^TK_j}M^{6T+4\sum\limits_{j=0}^TK_j}\big(1+\lambda M^2\big)^{2\sum\limits_{j=1}^TL_j+2\sum\limits_{j=0}^TK_j}}{\lambda^{\sum\limits_{j=1}^TL_j-T}\big(E\rho^2\big)^{\sum\limits_{j=1}^TL_j+2T
+2\sum\limits_{j=0}^TK_j}}\Bigg]<+\infty
\end{equation}
when $\lambda=\frac{\gamma}{dE\rho^2}$.

Then according to Lemma \ref{Lemma 4.1},
\[
\lambda_c(d)\leq \frac{\gamma}{dE\rho^2}
\]
for sufficiently large $d$ and hence
\[
\limsup_{d\rightarrow+\infty}d\lambda_c(d)\leq \frac{\gamma}{E\rho^2}
\]
for any $\gamma>1$.

Let $\gamma$ decrease to $1$, then we accomplish the proof.

\qed

\quad

Since we have shown that $\liminf_{d\rightarrow+\infty}d\lambda_c(d)\geq \frac{1}{E\rho^2}$ in Section \ref{section three}, the whole proof of Theorem \ref{theorem main 2.1} is completed.

\quad

\textbf{Acknowledgments.} The author is grateful to the financial support from the National Natural Science Foundation of China with grant number 11171342.

{}
\end{document}